\documentclass{article}
\usepackage{graphicx}
\usepackage{amsmath}
\usepackage{amssymb}
\title{\bfseries Laguerre polynomials and
 the inverse Laplace transform using discrete data}         
\date{} 

\begin{document}

 \maketitle
\author{\centerline{Tran Ngoc Lien\footnote{Faculty of Sciences,
Cantho University, Cantho, Vietnam. Email: tnlien@ctu.edu.vn.} ,
Dang Duc Trong\footnote{Department of Mathematics and Computer
Sciences, Hochiminh City National University, 227 Nguyen Van Cu,
Hochiminh City, Vietnam. Email: ddtrong@mathdep.hcmuns.edu.vn.  }
 and Alain Pham Ngoc Dinh\footnote{Mapmo UMR 6628, Universit\'e d'Orleans, BP 6759, 45067 Orleans Cedex, France. Email: alain.pham@math.cnrs.fr}  }     
 \ \par{\footnotesize{ {\bf Abstract.} We consider the problem of
finding a function defined on $(0,\infty)$ from a countable set of
values of its Laplace transform. The problem is severely ill-posed.
We shall use the expansion of the function in a series of Laguerre
polynomials to convert the problem in an analytic interpolation
problem. Then, using the coefficients of Lagrange polynomials we
shall construct a stable approximation solution. Error estimate is
given. Numerical results are produced. \vglue 0.2cm {\bf Keywords
and phrases.} inverse Laplace transform, Laguerre polynomials,
Lagrange polynomials, ill-posed problem, regularization. \vglue
0.2cm {\bf Mathematics Subject Classification 2000.} 44A10,
30E05,33C45. } }
\vglue 0.4cm

\noindent{\Large\bf 1. Introduction} \vglue 0.4cm Let
$L^2_\rho(0,\infty)$ be the space of real Lebesgue measurable functions
defined on $(0,\infty)$ such that
$$   \Vert f\Vert_{L^2_\rho}^2\equiv
\int_0^\infty | f(x)|^2e^{-x}dx<\infty. $$
This is a Hilbert space corresponding to the inner product
$$<f,g>=\int_0^\infty f(x)g(x)e^{-x}dx.$$

We consider the problem
of recovering a function $f\in L^2_\rho(0,\infty)$ satisfying the
equations
$$ {\cal L}f(p_j)\equiv  \int_0^\infty e^{-p_jx}f(x)dx=\mu_j  \eqno{(DIL)}$$
where $p_j\in (0,\infty), j=1,2,3,...$

Generally, we have  the classical problem of finding a function $f(x)$ from
its given image $g(p)$  satisfying
\begin{equation}
{\cal L}f(p)\equiv\int_0^\infty e^{-px}f(x)dx=g(p),
\label{1In}
\end{equation}
where $p$ is in a subset $\omega$ of the complex plane.
 We note that ${\cal L}f(p)$ is usually
an analytic function on a half plane $\{Re\ p>\alpha\}$
for an appropriate real number $\alpha$. Frequently, the image of a Laplace transform is known only
on a subset $\omega$ of the right half plane $\{Re\ p>\alpha\}$. Depending on the set $\omega$, we shall have
appropriate methods to construct the function $f$ from the values in the set
$$  \{{\cal L}f(p):\ p\in\omega\}. $$
Hence,  there
are no universal methods of inversion of the Laplace transform.

If the data $g(p)$ is given as a function
on a line
$(-i\infty+a,+i\infty+a)$ (i.e., $\omega=\{p: \ p=a+iy,
\ y\in {\bf R}\}$)
on the complex plane then
we can use the Bromwich inversion formula
([26], p. 67) to find the function $f(x)$.

If $\omega\subset\{ p \in {\bf R}:\ p>0\}$ then we have the problem
of real inverse Laplace transform. The right hand side is known only
on $(0,\infty)$ or a subset of $(0,\infty)$. In this case, the use
of the Bromwich formula is therefore not feasible. The literature on
the subject is impressed in both theoretical  and computational
aspects (see, e.g. [2,3,10,16,18,22]). In fact, if the data $g(p)$
is given exactly then, by the analyticity of $g$, we have many
inversion formulas (see,e.g., [3,7,8,20,21,23]). In [3], the author
approximate the function $f$ by
$$  f(t)\cong \sum_{k=0}^N b_k(a)d^k(e^xg(e^x))/dx^k$$
where $b_k(a)$ are calculated and tabulated regularization coefficients and $g$ is the given Laplace transform of $f$.
Another method is developped
by Saitoh and his group ([4,5,20,21]), where the function $f$ is approximated
by  integrals having the form
$$    u_N(t)=\int_0^\infty g(s)e^{-st}P_N(st)ds, \ \ \ \ N=1,2,...$$
where $P_N$ is known (see [5]). Using the Saitoh formula, we can
get directly error estimates.

However,
in the case of unexact data, we have a severely trouble by the ill-posedness of the problem.
In fact, a solution corresponding to the unexact
 data do not exist if the data is nonsmooth, and in the case of existence, these do not depend continuously on the given data (that are represented by the right hand side of the equalities). Hence, a regularization method is in order.  In [7], the authors  used the Tikhonov method
to regularize the problem. In fact, in this method, we can approximate
$u_0$ by functions $u_\beta$ satisfying
$$   \beta u_\beta + {\cal L}^*{\cal L}u_\beta=
{\cal L}^*g, \ \ \ \ \ \ \beta>0.$$
Since $\cal L$ is self-adjoint (cf. [7]), the latter equation
can be written as
$$  \beta u_\beta+\int_0^\infty \frac{u_\beta(s)}{s+t}ds
=\int_0^\infty e^{-st}g(s)ds.$$
The latter problem is well-posed.

 Although the inverse Laplace transform has a rich
literature, the papers devoted to the problem with
discrete data are scarce. In fact, from the analyticity
of ${\cal L}f(p)$,if
${\cal L}f(p)$  is known on a countable subset
of $\omega \subset\{Re\ p>\alpha\}$ accumulating at a point then ${\cal L}f(p)$
  is known on the whole
$\{Re\ p>\alpha\}$. Hence, generally, a set of discrete data is
enough for constructing an approximation function of $f$. It is a
moment problem. In [15], the authors presented some theorems on the
stabilization of the inverse Laplace transform. The Laplace image is
measured at  $N$ points to within some error $\epsilon$. This is
achieved by proving parallel stabilization results for a related
Hausdorff moment problem. For a construction of an approximate
solution of (DIL), we note that the sequence of functions
$(e^{-p_jx})$ is (algebraically) linear independent and moreover the
vector space generated by the latter sequence is dense in
$L^2(0,\infty)$. The method of truncated expansion as presented in
([6], Section 2.1) is applicable and we refer the reader to this
reference for full details. In [11, 13], the authors convert (DIL)
into a moment problem of finding a function in $L^2(0,1)$ and, then,
they use Muntz polynomials to construct an approximation for $f$.

Now, in the present paper, we shall convert (DIL) to an analytic
interpolation problem on the Hardy space of the unit disc. After
that, we shall use Laguerre polynomials and coefficients of Lagrange
polynomials to construct the function $f$. An approximation
corresponding to the non exact data and error estimate will be
given.

The remainder of the paper divided into two sections. In Section 2,
we convert our problem into an interpolation one and give a
uniqueness result. In Section 3, we shall give two regularization
results in the cases of exact data and non exact data. Numerical
comparisons with exact solution are given in the last section.\vglue
0.4cm

\noindent{\Large\bf 2. A uniqueness result} \vglue 0.4cm

 In this paper we shall use
Laguerre polynomials
$$   L_n(x)=\frac{e^x}{n!}\frac{d^n}{dx^n}(e^{-x}x^n).$$
We note that $\{L_n\}$ is a sequence of orthonormal polynomials on
$L^2_\rho(0,\infty)$. We note that (see [1], [9],
 page 67)
$$  \exp\left(\frac{xz}{z-1}\right) (1-z)^{-1}=
\sum_{n=0}^\infty L_n(x)z^n.$$
Hence, if we have the expansion
$$  f(x)=\sum_{n-0}^\infty a_nL_n(x)$$
then
$$\int_0^\infty f(x)\exp\left(\frac{xz}{z-1}\right) (1-z)^{-1}
e^{-x}dx= \sum_{n=0}^{\infty}a_nz^n.$$
It follows that
$$ \sum_{n=0}^{\infty}a_nz^n=
\int_0^\infty f(x)\exp\left(\frac{x}{z-1}\right) (1-z)^{-1}
dx.  $$
Put $\Phi f (z)=\sum_{n=0}^{\infty}a_nz^n,\ \alpha_j=1-1/p_j$, one has
$$   \Phi f(\alpha_j)=p_j\mu_j,$$
i.e., we have an interpolation problem of finding an
analytic function $\Phi f$ in the Hardy space $H^2(U)$. Here,
we denote by $U$ the unit disc of the complex plane
and  by $H^2(U)$  the Hardy space. In fact, we recall that $H^2(U)$ is the space of all functions $\phi$ analytic
in $U$ and if, $\phi\in H^2(U)$ has the expansion
$\phi(z)=\sum_{k=0}^\infty a_kz^k$ then
$$ \Vert \phi\Vert^2_{H^2(U)}=\sum_{k=0}^\infty \vert
a_k\vert^2
  =\frac{1}{2\pi}\int_0^{2\pi}|\phi(e^{i\theta})|^2d\theta. $$
 We can verify directly that the linear operator $\Phi$ is an isometry from $L^2_\rho$ onto
$H^2(U)$. In fact, we have
\vglue 0.4cm

{\bf Lemma 1} {\it
Let $f\in L^2_\rho(0,\infty)$. Then ${\cal L}f(z)$
is analytic on $\{z\in {\bf C}| \ Re z>1/2\}$.
If we  have an expansion
$$   f=\sum_{n=0}^\infty a_nL_n$$
then one has $\Phi f\in H^2(U)$ and
$$  \Vert \Phi f\Vert^2_{H^2(U)}=
\sum_{n=0}^\infty \vert a_n\vert^2= \Vert
f\Vert^2_{L^2_\rho(0,\infty)}.$$ Moreover, If we have in addition
that $\sqrt{x} f'\in L^2_\rho$ then
$$  \sum_{n=0}^\infty n\vert a_n\vert^2\leq
\Vert \sqrt{x} f'\Vert^2_{L^2_\rho}.$$
}
 \vglue 0.4cm
{\bf Proof} \vglue 0.2cm Putting $F_z(t)=e^{-zt}f(t)$, we have
$F_z\in L^2(0,\infty)$ for every $Re z>1/2$. Hence ${\cal
L}f(z)=\int_0^\infty F_z(t)dt$ is analytic for $Re z>1/2$. From the
definitions of $L^2_\rho(0,\infty)$ and $H^2(U)$, we have  the
isometry equality. Now we prove the second inequalities. We first
consider the case $f',f''$ in the space
$$  B=\{g \ Lebesgue\ measurable\ on\ (0,\infty)|\
\sqrt{x}g\in L^2_\rho(0,\infty)\}.$$
We have the expansion
$$  f=\sum_{n=0}^\infty a_nL_n$$
where $a_n=<f,L_n>$.

The function $y=L_n$ satisfies the following equation (see [17])
$$  xy''+(1-x)y'+ny=0$$
which gives
$$  (xe^{-x}y')'+nye^{-x}=0.  $$
It follows that
\begin{eqnarray*}
na_n &=& \int_0^\infty f(x)nL_n(x)e^{-x}dx\\
&=& -\int_0^\infty f(x)(xe^{-x}L'_n(x))'dx\\
&=&\int_0^\infty f'(x)xe^{-x}L'_n(x)dx\\
&=&-\int_0^\infty (f'(x)xe^{-x})'L_n(x)dx\\
&=&-\int_0^\infty (xf"(x)+f'(x)-xf'(x))L_n(x)e^{-x}dx\\
&=&-<xf"+f'-xf', L_n>.
\end{eqnarray*}
Since $L_n$ is an orthonormal basis, we have the Fourier expansion
$$  xf"+f'-xf'=\sum_{n=0}^\infty (-na_n)L_n.$$
Using the Parseval equality we have
$$<xf"+f'-xf', f>=\sum_{n=0}^\infty (-na_n)a_n.$$
It can be rewritten as
$$ \int_0^\infty (xe^{-x}f'(x))'f(x)dx
= -\sum_{n=0}^\infty na^2_n.$$
Integrating by parts, we get
$$  \int_0^\infty xe^{-x}\vert f'(x)\vert^2dx
= \sum_{n=0}^\infty na^2_n.$$ Now, for $f'\in B$ we choose $(f_k)$
such that $f'_k , f"_k\in B$ for every $k=1,2,...$ and $\sqrt{x} f'_k \,
(\textrm {resp.} f_k) \to \sqrt{x} f' \;(\textrm {resp.} f) $  in
$L^2_\rho$ as $k\to\infty$. Assume that
$$f_k= \sum_{n=0}^\infty a_{kn}L_n.$$
Then we have $$ \int_0^\infty xe^{-x}\vert f_k'(x)\vert^2dx =
\sum_{n=0}^\infty na^2_{kn}.$$
The latter equality involves for every N
\begin{equation}
 \sum_{n=0}^N na^2_{kn} \leq \Vert
\sqrt{x}f_k'\Vert^2_{L^2_\rho(0,\infty)}\label{2In}
\end{equation}

Since $ f_k  \to  f$ in $L^2_\rho$ as $k\to\infty $ we have that
$a_{kn } \to a_n $ as $k \to \infty,$ for each $n$. On the other
hand, we have $\sqrt{x} f'_k  \to \sqrt{x} f'$ in $L^2_\rho$ as
$k\to\infty .$ .  Therefore, letting $k\to \infty$ in (2) we get
$$ \sum_{n=0}^N na^2_{n} \leq \Vert \sqrt{x}f'\Vert^2_{L^2_\rho(0,\infty)}.$$
Letting $N\to\infty$ in the latter inequality, we get the desired inequality.
\hfill $\blacksquare$ \vglue 0.4cm

Using Lemma 1, one has a uniqueness result
\vglue 0.6cm
{\bf Theorem 1.} {\it
Let $p_j>1/2$ for every
$j=1,2,...$. If
$$  \sum_{p_j>1} \frac{1}{p_j}+ \sum_{1/2<p_j<1}
\frac{2p_j-1}{p_j}=\infty $$
then Problem (DIL) has at most one solution
in  $L^2_\rho(0,\infty).$}
\vglue 0.4cm
{\bf Proof}
\vglue 0.2cm
 Let
$f_1,f_2\in L^2_\rho(0,\infty)$ be two solutions of (DIL). Putting
$g=f_1-f_2$ then $g\in L^2_\rho(0,\infty)$ and ${\cal L}g(p_j)=0$.
It follows that $\Phi g(1-1/p_j)=0,\ j=1,2,...$ It follows that
$\alpha_j=1-1/p_j$ are zeros of $\Phi g$. We have $\Phi g\in H^2(U)$
and
$$\sum_{j=1}^\infty (1-|\alpha_j|)=
\sum_{p_j>1} \frac{1}{p_j}+ \sum_{1/2<p_j<1}
\frac{2p_j-1}{p_j}=\infty .$$ Hence we get $\Phi g\equiv 0$ (see,
e.g.,[19], page 308).
 It follows that $g\equiv 0.$
This completes the proof of Theorem 1. \hfill $\blacksquare$ \vglue
0.4cm \noindent{\Large\bf 3. Regularization and error estimates}
\vglue 0.4cm In the section, we assume that $(p_j)$ is a bounded
sequence, $p_j\not= p_k$ for every $j\not=k$. Without loss of generality, we shall assume that
$\rho=1$ is an accumulation  point of $p_j$ . In fact, if
$p_j$ has an accumulation point $\rho_0>1$ then, by putting ${\tilde
f}(x)=e^{-(\rho_0-1)x}f(x)$ and $p'_j=p_j-\rho_0+1$, we can
transform the problem to the one of finding ${\tilde f}\in
L^2_\rho(0,\infty)$ such that
$$   \int_0^\infty e^{-p'_jx}{\tilde f}(x)dx=\mu_j,
\ \ \ \ j=1,2,...   $$
in which  $p'_j$ has the accumulation point $\rho=1$. In fact, in Theorem 2 below, we
shall assume that $\left|1-\frac{1}{p_j}\right|\leq \sigma$ for every
$j=1,2,...$, where $\sigma$ is a given number.

We denote by $\ell^{(m)}_k(\nu)$ the coefficient of $z^k$
in the expansion of the Lagrange polynomial
$L_m(\nu)\ (\nu=(\nu_1,...,\nu_m))$ of degree
(at most) $m-1$ satisfying
$$   L_m(\nu)(z_k)=\nu_k,\ \ \ \ \ 1\leq k\leq m, $$
where $z_k=\alpha_k.$ If $\phi$ is an analytic function on $U$, we also denote
 $$L_m(\phi)=L_m(\phi(z_1),...,\phi(z_m)).$$
We define
$$  L_m^\theta (\nu)(z)=
\sum_{0\leq k\leq \theta(m-1)}\ell^{(m)}_k(\nu) z^k.$$
The polynomial $L_m^\theta (\nu)$ is called
a truncated Lagrange polynomial (see also [25]).
For every $g\in L^2_\rho(0,\infty)$, we put
\begin{eqnarray*}
T_ng&=&(p_1{\cal L}g(p_1),...,
p_n{\cal L} g(p_n)),\\
Tg&=&(p_1{\cal L}g(p_1),...,
p_n{\cal L} g(p_n),...)\in \ell^\infty.
\end{eqnarray*}
Here, we recall that $\alpha_n=1-1/p_n$. We
shall approximate the function $f$ by
$$   F_m=\Phi^{-1} L_m^\theta (T_mf) =\sum_{0\leq k\leq \theta(m-1)}
\ell^{(m)}_k(T_mf)L_k.$$
We shall prove that $F_m$ is an approximation of $f$. Before stating and proving
the main results, some remarks are in order.

We first recall the concept of regularization.
Let $f$ be an exact solution of (DIL), we recall that a sequence of
linear  operator $A_n: \ell^\infty\rightarrow L^2_\rho(0,\infty)$ is a
regularization sequence
(or a regularizer) of Problem (DIL)
if $(A_n)$ satisfies two following conditions (see, e.g., [14], page 25)

(R1) For each $n$, $A_n$ is bounded,

(R2) $\lim_{n\to\infty} \Vert A_n(Tf)-f\Vert=0.$

The number "n" is called the {\it regularization parameter}. As a consequence
of (R1), (R2), we can get

(R3) For $\epsilon>0$, there exists the functions
$n(\epsilon)$ and $\delta(\epsilon)$
such that $\lim_{\epsilon\to 0}n(\epsilon)=\infty$,
$\lim_{\epsilon\to 0}\delta(\epsilon)=0$ and that
$$   \Vert A_{n(\epsilon)}(\mu)-f\Vert\leq \delta(\epsilon)$$
for every $\mu\in \ell^\infty$ such that
$\Vert\mu-Tf\Vert_\infty<\epsilon$.

In the present paper, the operator $A_n$ is $\Phi^{-1} L_m^\theta$.
The number $\epsilon$ is the error between the exact
data $Tf$ and the measured data $\mu$. For a given error $\epsilon$, there
are infinitely many ways of choosing the regularization parameter $n(\epsilon)$.
In the present paper, we give an explicit form of $n(\epsilon)$.

Next, in our paper, we have the interpolation problem of reconstruction
the analytic function $\phi=\Phi f\in H^2(U)$ from a sequence of its values
$(\phi(\alpha_n))$. As known,
the convergence of
$L_m(\phi)$ to $\phi$ depends heavily on the properties
of the points $(\alpha_n)$. The Kalm\'ar-Walsh theorem
(see, e.g.,[12], page 65) shows that $L_m(\phi)\to \phi $ for every $\phi$
in $C(\overline{U})$ for all $\phi$ analytic in a {\it neighborhood} of $\overline{U}$
if and only if $(\alpha_n)$ is uniformly distributed in $\overline{U}$, i.e.,
$$\lim_{m\to\infty}\sqrt[m]{\max_{|z|\leq 1}|(z-\alpha_1)...(z-\alpha_m)|}=1.$$
The Fejer points and the Fekete points are the sequences of points satisfying the latter condition (see [12], page 67). The Kalm\'ar-Walsh fails if $C(\overline{U})$ is
replaced by $H^2(U)$ (see [25] for a counterexample). Hence, the Lagrange
polynomial cannot use to reconstruct $\phi$. In [12], we proved a theorem similar
to the  Kalm\'ar-Walsh theorem for the case of $H^2(U)$. In fact,
 the Lagrange polynomials
will convergence if we "cut off" some terms of the Lagrange polynomial. Especially,
in [12] and the present paper,
the points $(\alpha_n)$ are, in general, not uniformly distributed.

In Theorem 2, we shall verify the condition (R2). More precisely, we have
\vglue 0.4cm

{\bf Theorem 2} {\it Let $\sigma\in (0,1/3)$, let  $f\in L^2_\rho(0,\infty)$
 and let
 $p_j>1/2$ for $j=1,2,...$ satisfy
$$   \left\vert 1-\frac{1}{p_j}\right\vert
\leq \sigma.$$
Put $\theta_0$ be the unique solution of the equation
(unknown x)
$$\frac{2\sigma^{1-x}}{1-\sigma}=1.$$
Then for $\theta\in (0,\theta_0)$, one has

$$\Vert f-F_m\Vert^2_{L^2_\rho}
\longrightarrow 0\ \ \ {\rm as}\ m\to\infty. $$
If, we assume in addition that  $\sqrt{x}f'\in L^2_\rho(0,\infty)$ then
$$  \Vert f-F_m\Vert^2_{L^2_\rho}
\leq (1+m\theta)^2 \Vert f\Vert^2_{L^2_\rho}\left(
\frac{2\sigma^{1-\theta}}{1-\sigma}\right)^{2m}
 + \frac{1}{m\theta}\Vert \sqrt{x}
 f'\Vert^2_{L^2_\rho(0,\infty)}.$$
}
\vglue 0.4cm
{\bf Proof}
\vglue 0.2cm
We have in view of Lemma 1
\begin{equation}
\Vert f-F_m\Vert^2_{L^2_\rho}=
\sum_{0\leq k\leq \theta(m-1)}\vert
\delta^{(m)}_k\vert^2 +\sum_{k>\theta(m-1)}
\vert a_k\vert^2
\label{1}
\end{equation}
where $\delta^{(m)}_k=a_k-\ell^{(m)}_k(T_mf)$.
We shall give an estimate for $\delta^{(m)}_k$. In fact,
we have
$$\Vert \Phi f-L_m(T_mf)\Vert^2_{H^2(U)}=
\sum_{k=0}^{m-1} \vert\delta^{(m)}_k\vert^2 + \sum_{k=m}^\infty
\vert a_k\vert^2. $$ On the other hand, the Hermite representation
(see, e.g. [12], page 59, [24]) gives
$$  \Phi f(z)-L_m(T_mf)(z)=\frac{1}{2\pi i}
\int_{\partial U} \frac{\omega_m(z)(\Phi f)(\zeta)d\zeta}{
\omega_m(\zeta)(\zeta-z)}$$ where
$\omega_m(z)=(z-\alpha_1)...(z-\alpha_m).$ Now, if we denote by
$\sigma^{(m)}_{-1}= \sigma^{(m)}_{-2}=...=0$ and
\begin{eqnarray*}
\sigma^{(m)}_0 &=& 1\\
\sigma^{(m)}_r &=& \sum_{1\leq j_1<...<j_r\leq m}
\alpha_{j_1}...\alpha_{j_r}\ \ \ (1\leq r\leq m),\\
\beta_s^{(m)} &=& \frac{1}{2\pi i} \int_{\partial U} \frac{ \Phi
f(\zeta)d\zeta}{ \zeta^{s+1}\omega_m(\zeta)}
\end{eqnarray*}
then we can write in view of the Hermite representation
$$ \Phi f(z)-L_m(T_mf)(z)=\sum_{k=0}^\infty
\left(\sum_{r=0}^k
(-1)^{r}\sigma^{(m)}_{m-r}\beta_{k-r}^{(m)}\right)z^k .$$ From the
latter representation, one gets
$$   \delta^{(m)}_k=\sum_{r=0}^k
(-1)^{r}\sigma^{(m)}_{m-r}\beta_{k-r}^{(m)}, \ \ \ 0\leq k\leq
m-1.$$ Now, by direct computation, one has
$$  |\beta_s^{(m)}|\leq \frac{1}{2\pi}\int_0^{2\pi}
\frac{|\Phi f(e^{i\theta})|}{|\omega_m(e^{i\theta})|}d\theta.$$
But one has
$$ |\omega_m(e^{i\theta})|\geq (|e^{i\theta}|-|\alpha_1|)...
(|e^{i\theta}|-|\alpha_m|)\geq (1-\sigma)^m.$$
Hence
\begin{eqnarray*}
\vert\beta_s^{(m)}\vert&\leq& \frac{1}{2\pi(1-\sigma)^m}\int_0^{2\pi}
|\Phi f(e^{i\theta})|d\theta  \\
&\leq& \Vert \Phi f\Vert_{H^2(U)}(1-\sigma)^{-m}.
\end{eqnarray*}
We also have
$$ \vert\sigma^{(m)}_{m-r} \vert\leq
\sigma^{m-r}C^r_m\leq \sigma^{m-k}2^m, $$ where
$C^k_m=\frac{m!}{k!(m-k)!}$. Hence, we have
$$ \vert \delta^{(m)}_k\vert \leq (1+m\theta)
\Vert f\Vert_{L^2_\rho}\left(
\frac{2\sigma^{1-\theta}}{1-\sigma}\right)^m.$$

From the latter inequality, one has in view of (3)
$$\Vert f-F_m\Vert^2_{L^2_\rho}
\leq (1+m\theta)^2 \Vert f\Vert^2_{L^2_\rho}\left(
\frac{2\sigma^{1-\theta}}{1-\sigma}\right)^{2m}
 +
\sum_{k\geq m\theta}^\infty \vert a_k\vert^2. $$ For $\theta\in (0,
\theta_0)$, one has
$$   0<\frac{2\sigma^{1-\theta}}{1-\sigma}
<\frac{2\sigma^{1-\theta_0}}{1-\sigma}=1.$$
Hence, we have
$$ \lim_{m\to\infty}\Vert f-F_m\Vert^2_{L^2_\rho} =0 $$
as desired, since on the one hand we have the comparison between an
exponential with base $b<1$ and a power function and in the other
hand the remain of a convergent series $\sum_{k=0}^\infty \vert
a_k\vert^2.$

Now if $\sqrt{x} f'\in L^2_\rho(0,\infty)$ then one has since
$\frac{k}{m\theta}>1$ and from Lemma 1
\begin{eqnarray*}
 \sum_{k>m\theta}^\infty \vert a_k\vert^2
\leq \frac{1}{m\theta}\sum_{k=0}^\infty k\vert a_k\vert^2 \leq
\frac{1}{m\theta}\Vert \sqrt{x} f'\Vert^2_{L^2_\rho}.
\end{eqnarray*}
This completes the proof of Theorem 2. \hfill $\blacksquare$ \
\par\vglue 0.6cm Now, we consider the case of non-exact data. In Theorem 3, we shall
consider the condition (R3) of the definition of the regularization. Put
$$  D_m=\max_{1\leq n\leq m}\left
(\max_{\vert z\vert\leq R}\left\vert
\frac{\omega_m(z)}{(z-\alpha_n)\omega'_m(\alpha_n)}\right\vert
\right).
$$ Let $\psi:\ [0,\infty)\rightarrow {\bf R}$ be an increasing
function satisfying
$$   \psi(m)\geq mD_m,\ \ \ \ m=1,2,...$$
and
$$   m(\epsilon)= [\psi^{-1}(\epsilon^{-3/4})]-1    $$
where $[x]$ is the greatest integer $\leq x$. \vglue 0.4cm
\noindent{\bf Theorem 3.} {\it
 Let $\sigma\in (0,1/3)$, let  $f, \sqrt{x}f'\in L^2_\rho(0,\infty)$
 and let
 $p_j>1/2$ for $j=1,2,...$ satisfy
$$   \left\vert 1-\frac{1}{p_j}\right\vert
\leq \sigma.$$
Put $\theta_0$ be the unique solution of the equation
(unknown x)
$$\frac{2\sigma^{1-x}}{1-\sigma}=1.$$
Let $\epsilon>0$  and let $(\mu_j^\epsilon)$ be a measured
data of $({\cal L}f(p_j))$ satisfying
$$   \sup_j | p_j({\cal L}f(p_j)-\mu_j^\epsilon)|
<\epsilon.$$
Then for $\theta\in (0,\theta_0)$, one has
$$
\Vert f-\Phi^{-1}L^\theta_{m(\epsilon)}(\nu^\epsilon)
\Vert^2_{L^2_\rho} \leq 2(1+m(\epsilon)\theta)^2 \Vert
f\Vert^2_{L^2_\rho}\left(
\frac{2\sigma^{1-\theta}}{1-\sigma}\right)^{2m(\epsilon)}
 + \frac{2}{m(\epsilon)\theta}\Vert \sqrt{x}
 f'\Vert^2_{L^2_\rho}\\
+2\epsilon^{1/2}.
$$}
where $\nu_j^\epsilon=p_j\mu_j^\epsilon$ for $j=1,2,...$
\vglue 0.2cm
{\bf Proof}

We note that
$$ L_m(T_mf)(z)-L_m(\nu^\epsilon)(z)=
\sum_{j=1}^m (p_j\mu_j-\nu^\epsilon_j)
\frac{\omega_m(z)}{(z-\alpha_j)\omega'_m(\alpha_j)}.
$$
It follows that
$$   \Vert L_m(T_mf)-L_m(\nu^\epsilon)\Vert_\infty
\leq \epsilon mD_m.$$ Hence
$$  \Vert L^\theta_m(T_mf)-L^\theta_m(\nu^\epsilon)\Vert_{H^2(U)}\leq
\Vert L_m(T_mf)-L_m(\nu^\epsilon)\Vert_\infty \leq \epsilon mD_m.
$$ It follows by the isometry property of $\Phi$
\begin{eqnarray*}
\Vert f-\Phi^{-1}L^\theta_m(\nu^\epsilon)
\Vert^2_{L^2_\rho}
&\leq& 2\Vert f-F_m\Vert^2_{L^2_\rho}+
 2\Vert \Phi^{-1}L^\theta_m(T_mf)-
\Phi^{-1}L^\theta_m(\nu^\epsilon)
\Vert^2_{L^2_\rho}\\
&\leq& 2(1+m\theta)^2 \Vert f\Vert^2_{L^2_\rho}\left(
\frac{2\sigma^{1-\theta}}{1-\sigma}\right)^{2m}
 + \frac{2}{m\theta}\Vert \sqrt{x}
 f'\Vert^2_{L^2_\rho}\\
& &+2\epsilon^2m^2D_m^2.
\end{eqnarray*}
By choosing $m=m(\epsilon)$ we get the desired result. \hfill
$\blacksquare$  \ \par \vglue 0.4cm

\noindent{\Large\bf 4. Numerical results} \vglue 0.4cm We present
some results of numerical comparison between the function $ f(x)$
given in $ L^{2}_{p}(0,\infty  )$ and its approximated form $ F_{m}$
as it is stated in Theorem 2.\ \par \ \ \ First consider the
function $ f(x)={\rm e}^{-x}$ and its expansion in Laguerre series\
\par
\begin{equation}
{\rm e}^{-x}={\sum _{n\geq  0}  ^{}  {{{1}\over{2^{n+1}  }}
L_{n}(x)}} .\label{lag-num2}
\end{equation}
So in the Hardy space $ H^{2}(U)$, we have to interpolate the
analytic function\ \par
\begin{equation}
\Phi f(x)={\sum _{n\geq  0}  ^{}  {{{1}\over{2^{n+1}  }}
x^{n}={{1}\over{2-x}} }}  \label{lag-num3}
\end{equation}
by the Lagrange polynomial $ L_{m}(T_{m}f)$, interpolation defined
by\ \par
\begin{equation}
L_{m}(T_{m}f){\left( 1-{{1}\over{p_{i}  }}  \right) }  =p_{i}  {\int
_{0} ^{\infty  }  {{\rm e}^{\scriptstyle -p_{i}x}  {\rm
e}^{-x}dx={{p_{i} }\over{p_{i}+1}}  }}  \label{lag-num4}
\end{equation}
where $ p_{i}\longrightarrow  1$ as $ i\to \infty  $.\
\par On the interval $ (-1.8,+1.8)$ we have drawn in Fig.1 the
curves $ {\rm e}^{-x}$ and its approximation $ L_{m}(T_{m}f)(x)$ for
$ m=10$. If $ m=12$ there is divergence for our interpolation
(Fig.2) outside the interval $ (-1,+1).$

\par
\centerline{
\includegraphics[width=3in]{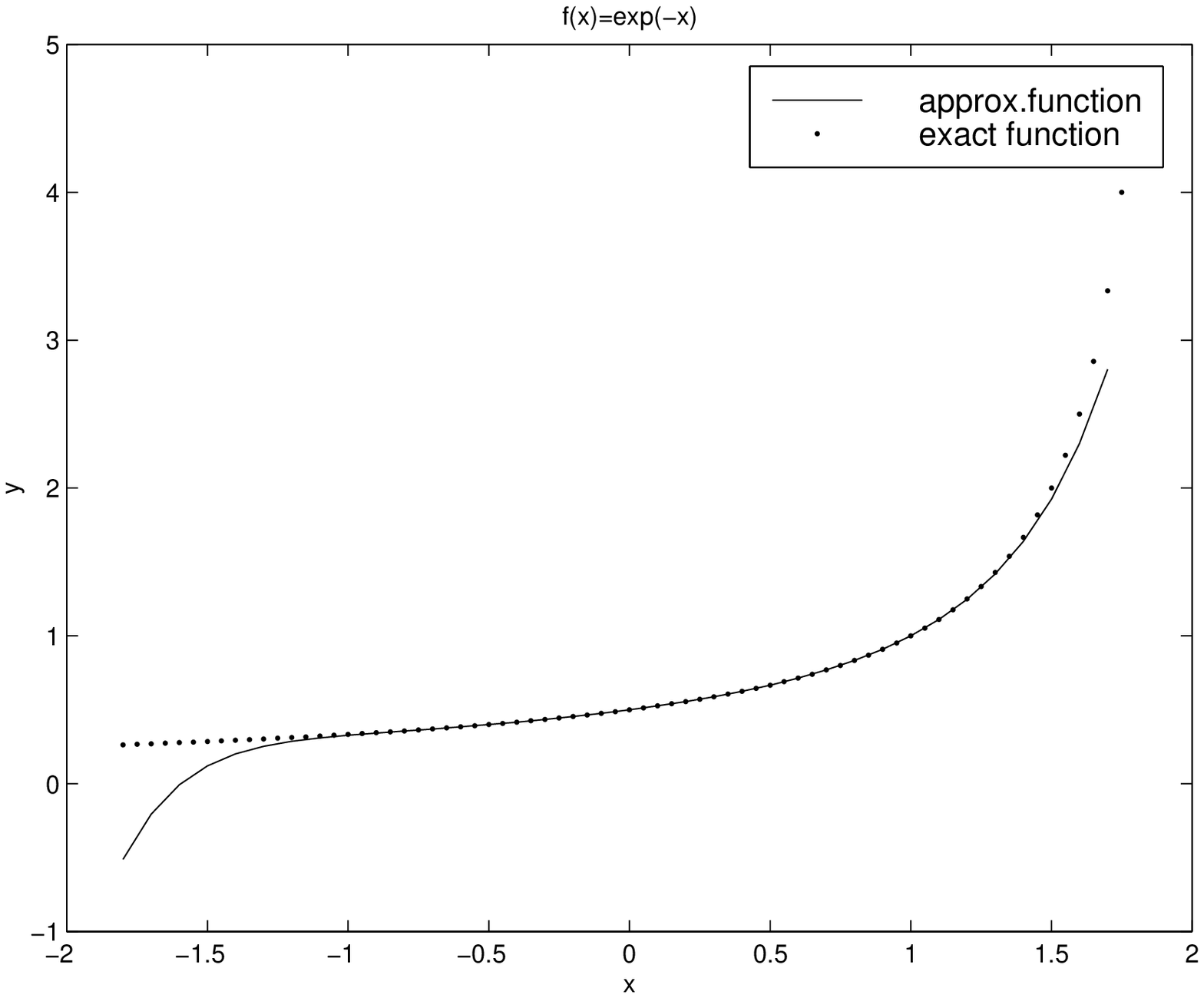}}
\centerline{Fig. 1}\par
\smallskip
\par
\centerline{
\includegraphics[width=3in]{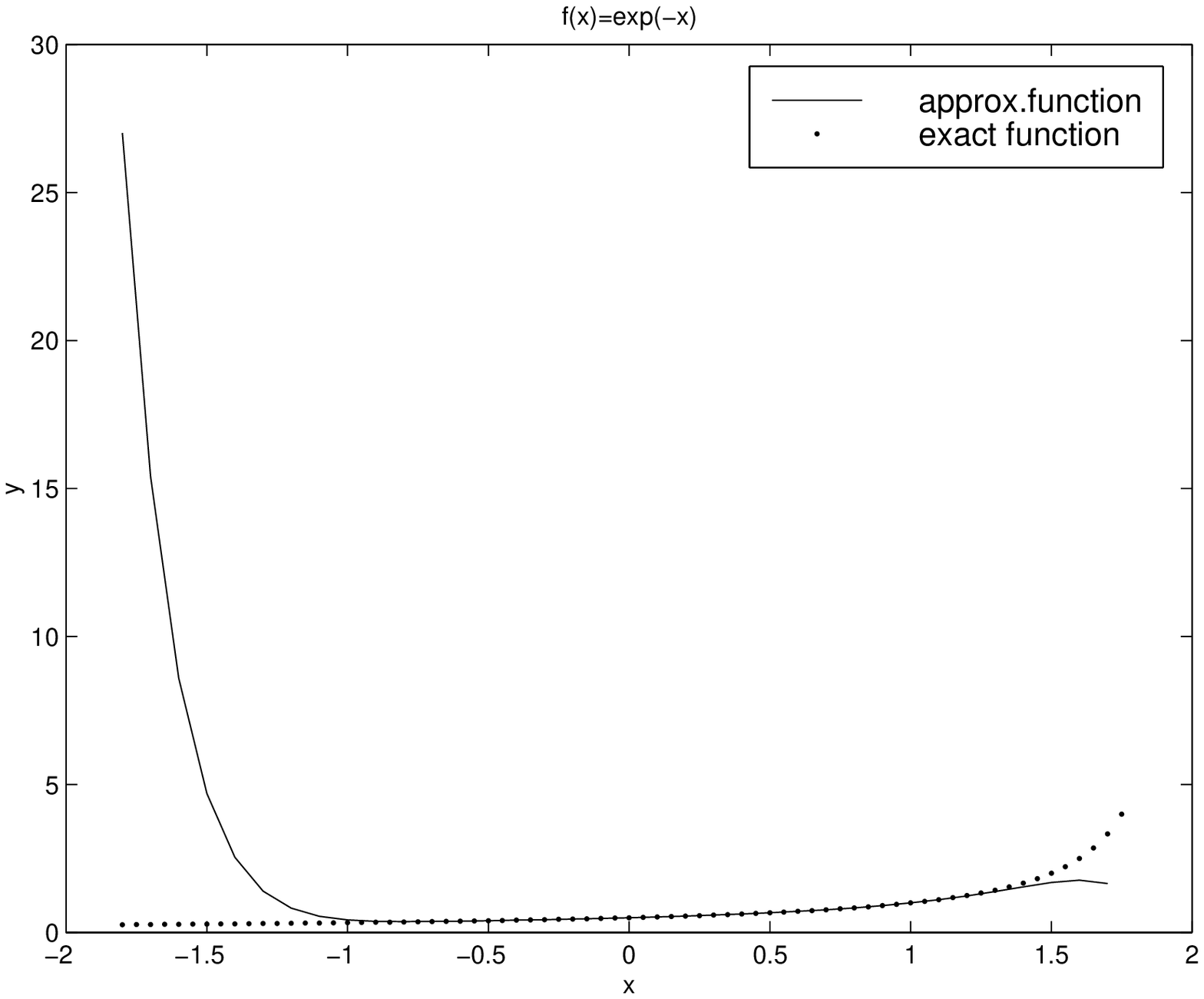}}
\centerline{Fig. 2}\par
\smallskip

\ \par In our 2nd example we have chosen the function\ \par
\begin{equation}
f(x)={\rm e}^{x/4}={{4}\over{3}}  {\sum _{n\geq  0}  ^{}  {{\left(
{{-1}\over{3}} \right) }  ^{\scriptstyle n}  L_{n}(x)}}
.\label{lag-num5}
\end{equation}
In the Hardy space the function \ \par
$$ \Phi f(x)={{4}\over{3}}  {\sum _{n\geq  0}  ^{}  {{\left( {{-x}\over{3}}
\right) }  ^{\scriptstyle n}  ={{4}\over{3+x}}  }}  $$ is
approximated by the Lagrange polynomial $ L_{m}(T_{m}f)$ at the
points $ (1-{{\displaystyle 1}\over{\displaystyle p_{i}}}
,{{\displaystyle -4p_{i}}\over{\displaystyle 1-4p_{i}}}  )$, $
p_{i}\longrightarrow  1$ as $ i\to  \infty $.\ \par The Fig.3
(resp.Fig.4) show the quite good convergence (resp.divergence) on
the interval $ (-2.8,2.8)$ with $ m=4$ (resp. $ m=11$).
\par
\centerline{
\includegraphics[width=3in]{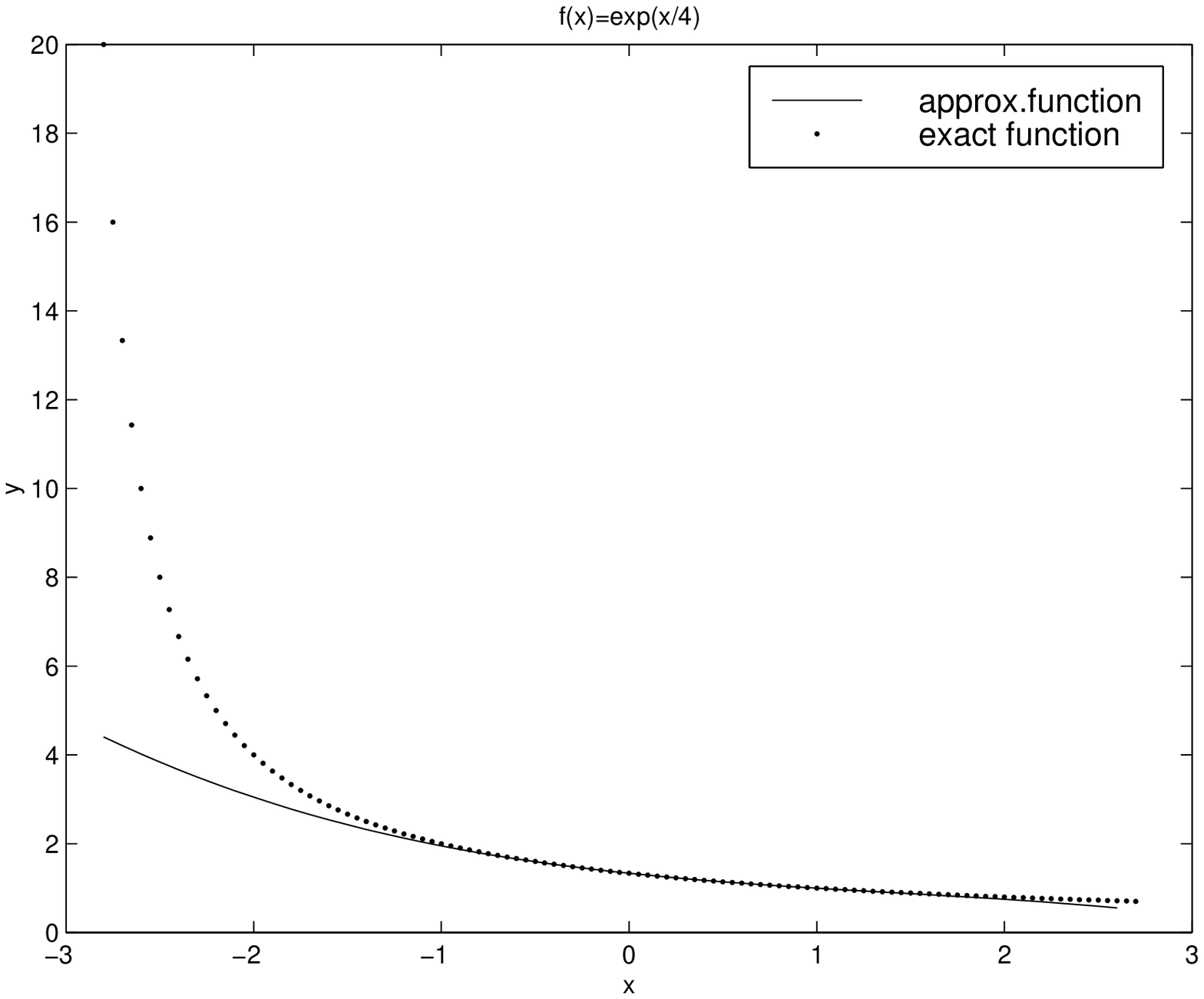}}
\centerline{Fig.3 }\par
\smallskip
\par
\centerline{
\includegraphics[width=3in]{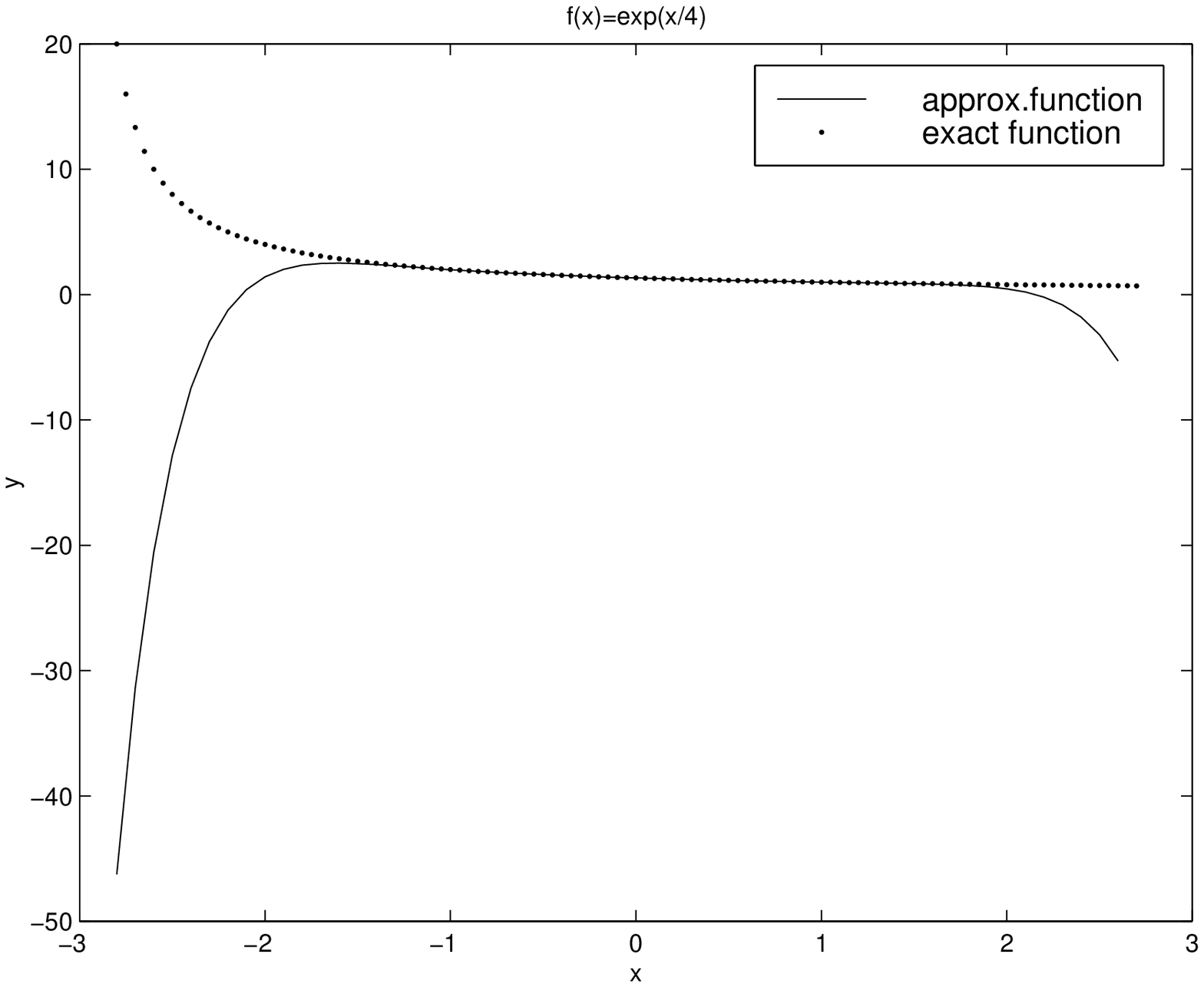}}
\centerline{Fig. 4}\par
\smallskip

In both cases we have chosen $ \theta  _{0}=0.29$ with $ \sigma
=0.25$ ($ \theta _{0}$ given by $ {{\displaystyle 2\sigma ^{1-\theta
o}}\over{\displaystyle 1-\sigma }}  =1,\  0<\sigma <{{\displaystyle
1}\over{\displaystyle 3}}  $). So in the 2nd case the truncated
Lagrange polynomial is almost verified since $ 11\times  0.29\sim
3.2$.\ \par

{\bf Acknowledgements. } The authors wish to thank the referees for
their pertinent remarks, leading to improvements in the original
manuscript.
\end{document}